\documentclass[reqno]{amsart}
\title[Legendre theorems for overpartitions]{Legendre theorems for certain overpartitions and overpartition pairs}
\usepackage{amssymb,amsmath,amsthm,epsfig,graphics,latexsym,hyperref}
\theoremstyle{definition}
\newtheorem{definition}{Definition}
\theoremstyle{plain}
\newtheorem{lemma}      {Lemma}

\newtheorem{theorem}    {Theorem}

\newtheorem{corollary}  {Corollary}

\theoremstyle{remark}

\numberwithin{equation}{section}

\newcommand{\fr}{\frac}

 \newmuskip\pFqskip
\mathchardef\pFcomma=\mathcode`, % keep a copy of the comma
\newcommand*\pFq[5]{%
  \begingroup
  \begingroup\lccode`~=`,
    \lowercase{\endgroup\def~}{\pFcomma\mkern\pFqskip}%
  \mathcode`,=\string"8000
  {}_{#1}\phi_{#2}\biggl[\genfrac..{0pt}{}{#3}{#4};#5\biggr]%
  \endgroup
}
\begin{document}
\author[ G. E. Andrews and M. El Bachraoui]{George E. Andrews and Mohamed El Bachraoui}
\address{The Pennsylvania State University, University Park, Pennsylvania 16802}
\email{andrews@math.psu.edu}
\address{Dept. Math. Sci,
United Arab Emirates University, PO Box 15551, Al-Ain, UAE}
\email{melbachraoui@uaeu.ac.ae}
%American University of Sharjah, P.O. Box 26666, Sharjah, UAE}
%\email{jgriffin@aus.edu}
%
\begin{abstract}
Motivated by two Legendre-type formulas for overpartitions, we derive a variety of their companions as Legendre theorems for overpartition pairs.
This leads to equalities of subclasses of overpartitions and overpartition pairs.
\end{abstract}
\date{\textit{\today}}
\thanks{First author partially supported by Simons Foundation Grant 633284}
\maketitle
\section{Introduction}\label{sec introduction}
Throughout $q$ denotes a complex number satisfying $|q|<1$, $m$ and $n$ denote nonnegative integers.
%$\mathbb{N}$ is the set of positive integers,
%$\mathbb{N}_0=\mathbb{N}\cup\{0\}$, and $\mathbb{Z}$ is the set of integers.
We will use the following standard notation for $q$-series~\cite{Andrews, Gasper-Rahman}
\[
(a;q)_0 = 1,\  (a;q)_n = \prod_{j=0}^{n-1} (1-aq^j),\quad
(a;q)_{\infty} = \prod_{j=0}^{\infty} (1-aq^j),
\]
\[
(a_1,\ldots,a_k;q)_n = \prod_{j=1}^k (a_j;q)_n,\ \text{and\ }
(a_1,\ldots,a_k;q)_{\infty} = \prod_{j=1}^k (a_j;q)_{\infty}.
\]
We will frequently use without reference the following basic facts of $q$-series~\cite{Andrews, Gasper-Rahman}
\begin{equation}\label{basic-facts}
(a;q)_{n+m} = (a;q)_{m} (aq^{m};q)_n,\
(a;q)_{\infty} = (a;q)_n (aq^n;q)_{\infty},\
(a;q)_{\infty} = (a;q^2)_{\infty}(aq;q^2)_{\infty}.
\end{equation}

Letting $p_e(\mathcal{D},n)$ (resp. $p_o(\mathcal{D},n)$) denote the number of partitions of $n$ into an
even (resp. odd) number of distinct parts, it is easy to see that their difference has the following generating function~\cite{Andrews}
\begin{equation}\label{Legendre-0}%~\cite{Andrews}
\sum_{n\geq 0} \big( p_e(\mathcal{D},n) - p_o(\mathcal{D},n) \big) q^n =(q;q)_\infty.
\end{equation}
Then combining~\eqref{Legendre-0} with Euler's pentagonal theorem~\cite{Andrews}
\begin{equation}\label{Euler-1}
(q;q)_\infty = \sum_{n=-\infty}^\infty (-1)^n q^{\fr{3n^2-n}{2}}
=1+\sum_{n\geq 1}(-1)^n q^{\fr{3n^2+n}{2}} +\sum_{n\geq 1}(-1)^n q^{\fr{3n^2-n}{2}},
\end{equation}
we get
\begin{equation}\label{Legendre-1}
\sum_{n\geq 0} \big( p_e(\mathcal{D},n) - p_o(\mathcal{D},n) \big) q^n
=\sum_{n=-\infty}^\infty (-1)^n q^{\fr{3n^2-n}{2}}
\end{equation}
which is equivalent to Legendre's~\cite{Legendre} celebrated result
\begin{equation}\label{Legendre}
p_e(\mathcal{D},n) - p_o(\mathcal{D},n)
=\begin{cases} (-1)^k, & \text{if\ }n=k(3k\pm 1)/2\ \text{for some\ } k\in\mathbb{Z}, \\
0, & \text{otherwise.}
\end{cases}
\end{equation}
Formulas of type~\eqref{Legendre-1} or~\eqref{Legendre} are referred to as Legendre theorems.
An immediate consequence of~\eqref{Legendre} is the vanishing of the difference
$p_e(\mathcal{D},n) - p_o(\mathcal{D},n)$ and obviously the equality of $p_e(\mathcal{D},n)$ and $p_o(\mathcal{D},n)$ at any positive integer
$n$ which is not a pentagonal number.

An overpartition~\cite{Corteel-Lovejoy} of $n$ is a partition of $n$ where the first
occurrence of each part may be overlined. The number of overpartitions of $n$, written
$\overline{p}(n)$, has the following generating function
\[
\sum_{n=0}^\infty \overline{p}(n) q^n = \fr{(-q;q)_\infty}{(q;q)_\infty}.
\]
Note that overlined parts in overpartitions are distinct by definition.
We say that an overpartition has distinct parts if its non-overlined parts are distinct too.
Letting
$\overline{p}_d(n)$ denote the number of overpartitions of $n$ into distinct parts, it is easy to see that
\[
\sum_{n=0}^\infty \overline{p}_d(n) q^n = (-q;q)_\infty^2.
\]
Recently, The first author and A. J. Yee~\cite{Andrews-Yee} did an extensive study of Legendre theorems for overpartitions revealing that this was a topic filled with elegant possibilities.
For instance, letting $TH(n)$ denote the number of overpartitions of $n$ in which there is both
an overlined and a non-overlined largest part and letting $THE(n)$ denote the number of $TH$-overpartitions
with an even number of parts minus the number with an odd number of parts, the authors proved that
\begin{equation}\label{AndYee id}
THE(n)
=\begin{cases}
(-1)^n(2k-1),\ \text{if\ } k^2<n< (k+1)^2\ \text{for some\ } k\in\mathbb{N}, \\
(-1)^n(2k-2),\ \text{if\ } n=k^2\ \text{for some\ } k\in\mathbb{N}.
\end{cases}
\end{equation}
The second author in~\cite{Bachraoui 2023} studied some overpartitions into distinct parts and
obtained Legendre theorems analogue
to~\eqref{AndYee id} with the conditions involving squares $k^2$ replaced by conditions involving triangular numbers $\binom{k}{2}$.

An overpartition pair~\cite{Lovejoy 2006} of $n$ is a pair of overpartitions $\pi=(\lambda_1, \lambda_2)$ where the sum
of all of the parts is $n$. The number of all overpartition pairs of $n$, written $\overline{pp}(n)$, has the following
generating function
\[
\sum_{n=0}^\infty \overline{pp}(n) q^n = \fr{(-q;q)_\infty^2}{(q;q)_\infty^2}.
\]
We say that an overpartition pair $(\lambda_1, \lambda_2)$ has distinct parts if both $\lambda_1$
and $\lambda_2$ have distinct parts.
While there is a rich literature dealing with results of Legendre type along with non-negativity results and inequalities
for both partitions and overpartitions (see for instance~\cite{Andrews 2013, Bachraoui 2023-b, Berkovich-Grizzell 2014, Berkovich-Uncu 2019, Kim-Kim-Lovejoy 2020, Kim-Kim-Lovejoy 2021, Lovejoy 2005}), there seems not much to have been done in this direction for overpartition pairs.
Besides, partition pairs have an established and beautiful theory.  W. H. Burge in a seminal paper~\cite{Burge} studied partition pairs at length and revealed truly surprising theorems and rich possibilities for further research.
In this paper, we hope to follow the natural direction suggested by joint consideration of~\cite{Burge},~\cite{Andrews-Yee}, and~\cite{Bachraoui 2023}
as we shall examine Legendre type theorems related to overpartitions and pairs of overpartitions.  Namely, we will prove the following main results which are particularly lovely examples of Legendre type theorems arising in this area.
\begin{theorem}\label{thm FG'}
There holds
\[
\begin{split}
(a)\ \sum_{n=1}^\infty F'(n) q^n &:= \sum_{n=1}^\infty q^{n}(q^{n+1};q)_\infty (q^{n};q)_{n} = \sum_{n\geq 1}(-1)^{n+1} q^{\fr{3n^2-n}{2}}, \\
(b)\ \sum_{n=1}^\infty G'(n) q^n &:= \sum_{n=1}^\infty q^{n}(q^{n+1};q)_\infty (-q^{n};q)_{n} = \sum_{n\geq 1}(-1)^{n+1} q^{\fr{3n^2-n}{2}} + 2\sum_{n\geq 0} q^{6n^2+7n+2}.
\end{split}
\]
\end{theorem}
\begin{theorem}\label{thm A'}
We have
\[
\sum_{n=1}^\infty A'(n) q^n
:=\sum_{n=1}^\infty q^{n}(q^{n+1};q)_\infty^3  (q^{n};q)_{n}
%= q - 2q^3+3q^6-4q^{10}+5 q^{15}-6q^{21}+\cdots
=\sum_{n=1}^\infty (-1)^{n+1} n  q^{\fr{n(n+1)} {2}}.
\]
\end{theorem}
\begin{theorem}\label{thm A''}
We have
\[
\sum_{n=1}^\infty A''(n) q^n
:=\sum_{n=1}^\infty q^{n}(-q^{n+1};q)_\infty^2 (q^{n+1};q)_\infty  (q^{n};q)_{n}
%= q + 2q^3+3q^6+4q^{10}+5 q^{15}+6q^{21}+\cdots
=\sum_{n=1}^\infty  n  q^{\fr{n(n+1)} {2}}.
\]
\end{theorem}
\begin{theorem}\label{thm B'}
We have
\[
\sum_{n=1}^\infty B'(n) q^n
:=\sum_{n=1}^\infty q^{n}(-q^{n+1};q)_\infty^2 (q^{n+1};q)_\infty  (q^{n+1};q)_{n}
%= \Big(\fr{(q^2;q^2)_\infty}{(q;q^2)_\infty} \Big)^2 - \fr{(q^2;q^2)_\infty}{(q;q^2)_\infty}
= \Big(\sum_{n=0}^\infty q^{\fr{n(n+1)}{2}} \Big)^2 - \sum_{n=0}^\infty q^{\fr{n(n+1)}{2}}.
\]
\end{theorem}
Our results include the following two companions of the above identities.
\begin{theorem}\label{thm C'}
Let
\[
\sum_{n=1}^\infty C'(n) q^n
:=\sum_{n=1}^\infty q^{n}(q^{n+1};q)_\infty^2 (q^{n};q)_\infty  (q^{n};q)_{n}.
\]
Then we have
\[
\sum_{n=1}^\infty C'(n) q^{8n+2}
= 2 \sum_{r,n=0}^\infty (-1)^{n+1} q^{(2r+3)^2 + (2r+3+2n)^2} + \sum_{n=0}^\infty (-1)^{n+1} q^{1+(2n+1)^2}
\]
\[
+ \sum_{n=1}^\infty q^{2(2n+1)^2}
\]
\end{theorem}
\begin{theorem}\label{thm D'}
Let
\[
\sum_{n=1}^\infty D'(n) q^n
:=\sum_{n=1}^\infty q^{2n}(q^{n+1};q)_\infty^3  (q^{n};q)_{n}.
\]
Then we have
\[
\sum_{n=1}^\infty D'(n) q^{8n+2}
= 2 \sum_{r,n=0}^\infty (-1)^{n} q^{(2r+1)^2 + (2r+1+2n)^2} - \sum_{n=0}^\infty (-1)^{n}(n+1) q^{1+(2n+1)^2}
\]
\[
- \sum_{n=0}^\infty q^{2(2n+1)^2}
\]
\end{theorem}
%The generating function for the partition function $p(n)$ is~\cite{Andrews}
%\[
%\fr{1}{(q;q)_\infty} = \sum_{n=0}^\infty p(n) q^n.
%\]
%
%It turns out that some of our difference formulas vanish at the non-triangular numbers and some of them vanish at the positive integers
%that are not sums of two triangular numbers. This leads to equalities of subclasses of overpartition pairs.
%Proofs for our results rely on manipulations with $q$-series
%and we do not know bijective proofs for them.
% Our proofs combine basic facts from the theory of hypergeometric $q$-series with the Bailey pair machinery.

The rest of the paper is organized as follows. In Section~\ref{sec combinatorial} we give combinatorial interpretations in terms of overpartitions
and pairs of overpartitions for our sequences.
Sections~\ref{sec proof FG'}-\ref{sec proof D'} are devoted to the proofs of
the main theorems.
Finally in Section~\ref{sec conclusion} we close by some remarks and questions suggested by this work.
\section{Combinatorial interpretations}\label{sec combinatorial}
Throughout all of the overpartitions and overpartition pairs have a smallest part. We will write $s(\pi)$ to denote
the smallest part of the partition $\pi$.
The two sequences $F'(n)$ and $G'(n)$ in Theorem~\ref{thm FG'} have the following natural interpretations as
overpartition differences.
\begin{definition}\label{def FG'}
For any positive integer $n$ let $F(n)$ denote the number of overpartitions $\pi$ of
$n$ into distinct parts where $s(\pi)$ occurs overlined and the non-overlined parts
are in the half-open interval $[s(\pi), 2 s(\pi))$.
Let $F_0(n)$ (resp. $F_1(n)$) denote the number of overpartition counted by $F(n)$ in which the number of
parts is even (resp. odd).
By letting the term $q^n (q^{n+1};q)_{\infty}$ generate the overlined parts of $\pi$
and $(q^{n};q)_{n} $ generate its non-overlined parts,
it is easy to check that
\begin{equation}\label{gen F'}
\sum_{n=1}^\infty \big(F_1(n)-F_0(n)\big) q^n
=\sum_{n=1}^\infty q^{n}(q^{n+1};q)_\infty (q^{n};q)_{n},
\end{equation}
which by the definition of $F'(n)$ in Theorem~\ref{thm FG'} yields
\[
F'(n) = F_1(n)-F_0(n).
\]
Similarly, let $G_0(n)$ (resp. $G_1(n)$) denote the number of overpartitions counted by $F(n)$ in which the number of overlined
parts is even (resp. odd).
By letting the term $q^n (q^{n+1};q)_{\infty}$ generate the overlined parts of $\pi$
and $(-q^{n};q)_{n} $ generate its non-overlined parts,
it is easy to check that
\begin{equation}\label{gen G'}
\sum_{n=1}^\infty \big(G_1(n)-G_0(n)\big) q^n
=\sum_{n=1}^\infty q^{n}(q^{n+1};q)_\infty (-q^{n};q)_{n}.
\end{equation}
That is,
\[
G'(n) = G_1(n)-G_0(n).
\]
\end{definition}
For example, $F(4)=4$ counting
\[
\bar{4}, \bar{3}+\bar{1}, \bar{2}+2, \bar{2}+\bar{1}+1.
\]
We have $F_0(4)=2$ counting $\bar{3}+\bar{1}$ and $\bar{2}+2$ and we have $F_1(4)=2$ counting $\bar{4}$ and $\bar{2}+\bar{1}+1$ and thus
$F'(4)=0$. Furthermore, it is easily checked that $G'(4)=0$. This agrees with Theorem~\ref{thm FG'}.
We now focus on the sequences $A'(n)$, $A''(n)$, $B'(n)$, $C'(n)$, and $D'(n)$
whose natural interpretations are differences of overpartition pairs.
\begin{definition}\label{def A}
For any positive integer $n$ let $A(n)$ denote the number of overpartition pairs $\pi=(\lambda_1, \lambda_2)$ of
$n$ into distinct parts where $s(\pi)=s(\lambda_1)$ and $s(\lambda_1)$ occurs overlined, the non-overlined parts of $\lambda_1$
are $>s(\pi)$ and the non-overlined parts of $\lambda_2$ are in the half-open interval $[s(\pi), 2 s(\pi))$.
Let $A_0(n)$ (resp. $A_1(n)$) denote the number of overpartition pairs counted by $A(n)$ in which the number of
parts is even (resp. odd) and let
\[
A'(n) = A_1(n)-A_0(n).
\]
By letting the term $q^n (q^{n+1};q)_{\infty}$ generate the overlined parts of $\lambda_1$
and $(q^{n+1};q)_{\infty} $ generate its non-overlined parts and letting
the term $(q^{n+1};q)_{\infty} $ generate the overlined parts of $\lambda_2$
and $(q^{n};q)_{n} $ generate its non-overlined parts,
it is easy to check that
\begin{equation}\label{gen A'}
\sum_{n=1}^\infty A'(n) q^n
=\sum_{n=1}^\infty q^{n}(q^{n+1};q)_\infty^3  (q^{n};q)_{n}.
\end{equation}
Furthermore, let $A_2(n)$ (resp. $A_3(n)$) denote the number of overpartition pairs counted by $A(n)$ in which the number of
non-overlined parts is even (resp. odd) and let
\[
A''(n) = A_2(n)-A_3(n).
\]
We have
\begin{equation}\label{gen A''}
\sum_{n=1}^\infty A''(n) q^n
=\sum_{n=1}^\infty q^{n}(-q^{n+1};q)_\infty^2 (q^{n+1};q)_\infty  (q^{n};q)_{n}.
\end{equation}
\end{definition}
For example, $A(3)=4$ with relevant partition pairs
\[
\big( (\bar{3}),\emptyset \big), \big( (\bar{2}, \bar{1}),\emptyset \big), \big( (2, \bar{1}),\emptyset \big),
\ \text{and\ } \big( (\bar{1}),(\bar{2}) \big).
\]
We have $A_0(3)=1$ with the only overpartition pair $\big( (\bar{3}),\emptyset \big)$
and $A_1(3)=3$ counting
\[
\big( (\bar{2}, \bar{1}),\emptyset \big), \big( (2, \bar{1}),\emptyset \big),
\ \text{and\ } \big( (\bar{1}),(\bar{2}) \big)
\]
and thus $A'(3)= 1-3=-2$.
On the other hand,
we have
$A_2(3)=3$ counting
\[
\big( (\bar{3}),\emptyset \big), \big( (\bar{2}, \bar{1}),\emptyset \big),
\ \text{and\ } \big( (\bar{1}),(\bar{2}) \big)
\]
and $A_3(3)=1$ with the only relevant pair $\big( (2, \bar{1}),\emptyset \big)$,
and thus $A''(3)= 3-1=2$.
\begin{definition}\label{def B}
For any positive integer $n$ let $B(n)$ denote the number of overpartition pairs $\pi=(\lambda_1, \lambda_2)$ of
$n$ into distinct parts where $s(\pi)=s(\lambda_1)$ is overlined and occurs exactly once in $\pi$ and the non-overlined parts in $\lambda_2$
are $\leq 2 s(\pi)$.
Let $B_0(n)$ (resp. $B_1(n)$) denote the number of overpartition pairs counted by $B(n)$ in which the number of
non-overlined parts is even (resp. odd) and let
\[
B'(n) = B_0(n)-B_1(n).
\]
By letting the term $q^n (-q^{n+1};q)_{\infty}$ generate the overlined parts of $\lambda_1$
and $(q^{n+1};q)_{\infty} $ generate its non-overlined parts and letting
the term $(-q^{n+1};q)_{\infty} $ generate the overlined parts of $\lambda_2$
and $(q^{n+1};q)_{n} $ generate its non-overlined parts,
it is directly verified that
\begin{equation}\label{gen B'}
\sum_{n=1}^\infty B'(n) q^n
=\sum_{n=1}^\infty q^{n}(-q^{n+1};q)_\infty^2 (q^{n+1};q)_\infty  (q^{n+1};q)_{n}.
\end{equation}
\end{definition}
For example, $B(3)=5$ counting
\[
\big( (\bar{3}),\emptyset \big), \big( (\bar{2},\bar{1}),\emptyset \big), \big( (2,\bar{1}),\emptyset \big),
\big( (\bar{1}),(\bar{2}) \big), \big( (\bar{1}),(2) \big).
\]
We have $B_0(3)=3$ counting
\[
\big( (\bar{3}),\emptyset \big), \big( (\bar{2},\bar{1}),\emptyset \big),
\big( (\bar{1}),(\bar{2}) \big),
\]
and $B_0(3)=2$ counting
\[
\big( (2,\bar{1}),\emptyset \big)\ \text{and\ }
\big( (\bar{1}),(2) \big)
\]
and thus $B'(3)= 3-2=1$.
\begin{definition}\label{def C}
For any positive integer $n$ let $C(n)$ denote the number of overpartition pairs $\pi=(\lambda_1, \lambda_2)$ of
$n$ into distinct parts where $s(\pi)=s(\lambda_1)$ occurs overlined, the overlined parts of $\lambda_2$
are $>s(\pi)$ and its non-overlined parts $< 2 s(\pi)$.
Let $C_0(n)$ (resp. $C_1(n)$) denote the number of overpartition pairs counted by $C(n)$ in which the number of
parts is even (resp. odd) and let
\[
C'(n) = C_1(n)-C_0(n).
\]
By letting the term $q^n (q^{n+1};q)_{\infty}$ generate the overlined parts of $\lambda_1$
and $(q^{n};q)_{\infty} $ generate its non-overlined parts and letting
the term $(q^{n+1};q)_{\infty} $ generate the overlined parts of $\lambda_2$
and $(q^{n};q)_{n} $ generate its non-overlined parts,
it is easily seen that
\begin{equation}\label{gen C'}
\sum_{n=1}^\infty C'(n) q^n
=\sum_{n=1}^\infty q^{n}(q^{n+1};q)_\infty^2 (q^{n};q)_\infty  (q^{n};q)_{n}.
\end{equation}
\end{definition}
For example, $C(3)=5$ with relevant partition pairs
\[
\big( (\bar{3}),\emptyset \big), \big( (\bar{2},\bar{1}),\emptyset \big), \big( (2,\bar{1}),\emptyset \big),
\ \text{and\ } \big( (\bar{2}),(\bar{1})\big), \big( (\bar{1},1),(1) \big).
\]
We have $C_0(3)=3$ counting
\[
\big( (\bar{2},\bar{1}),\emptyset \big), \big( (2,\bar{1}),\emptyset \big),
\ \text{and\ } \big( (\bar{2}),(\bar{1})\big)
\]
and $C_1(3)=2$ counting
\[
\big( (\bar{3}),\emptyset \big)\ \text{and\ } \big( (\bar{1},1),(1) \big),
\]
and thus $C'(3)= 2-3=-1$.
%
%Our last overpartition pairs are a slight modification of Definition~\ref{def A}.
Our last example of differences of overpartition pairs is a slight modification of Definition~\ref{def A}. In particular, the smallest must
appear in both components of the overpartition pair $\pi=(\lambda_1,\lambda_2)$, i.e. $s(\pi)=s(\lambda_1)=s(\lambda_2)$.
\begin{definition}\label{def D}
For any positive integer $n$ let $D(n)$ denote the number of overpartition pairs $\pi=(\lambda_1, \lambda_2)$ of
$n$ into distinct parts where $s(\pi)=s(\lambda_1)=s(\lambda_2)$ occurs overlined, the non-overlined parts of $\lambda_1$
are $>s(\pi)$ and the non-overlined parts of $\lambda_2$ are in the interval $[s(\pi), 2 s(\pi))$.
Let $D_0(n)$ (resp. $D_1(n)$) denote the number of overpartition pairs counted by $D(n)$ in which the number of
parts is even (resp. odd) and let
\[
D'(n) = D_0(n)-D_1(n).
\]
By letting the term $q^n (q^{n+1};q)_{\infty}$ generate the overlined parts of $\lambda_1$
and $(q^{n+1};q)_{\infty} $ generate its non-overlined parts and letting
the term $q^n(q^{n+1};q)_{\infty} $ generate the overlined parts of $\lambda_2$
and $(q^{n};q)_{n} $ generate its non-overlined parts,
it is easy to check that
\begin{equation}\label{gen D'}
\sum_{n=1}^\infty D'(n) q^n
=\sum_{n=1}^\infty q^{2n}(q^{n+1};q)_\infty^3  (q^{n};q)_{n}.
\end{equation}
\end{definition}
For example, $D(4)=5$ with relevant partition pairs
\[
\big( (\bar{2}),(\bar{2}) \big), \big( (\bar{2}, \bar{1}),(\bar{1}) \big), \big( (2, \bar{1}),(\bar{1}) \big),
\ \text{and\ } \big( (\bar{1}),(\bar{2},\bar{1})\big).
\]
We have $D_0(4)=1$ counting $\big( (\bar{2}),(\bar{2}) \big)$
and $D_1(4)=3$ counting
\[
\big( (\bar{2}, \bar{1}),(\bar{1}) \big), \big( (2, \bar{1}),(\bar{1}) \big),
\ \text{and\ } \big( (\bar{1}),(\bar{2},\bar{1})\big),
\]
and thus $D'(4)= 1-3=-2$.

Similarly, we have $D(5)=6$, enumerating
\[
\big( (\bar{3},\bar{1}),(\bar{1}) \big), \big( (3,\bar{1}),(\bar{1}) \big), \big( (\bar{1}),(\bar{3},\bar{1}) \big),
\big( (\bar{2},\bar{1}),(\bar{1},1) \big), \big( (2,\bar{1}),(\bar{1},1) \big),
\ \text{and\ } \big( (\bar{1}),(\bar{2},\bar{1},1) \big)
\]
and we can easily see that $D_0(5) = D_1(5)=3$ and thus $D'(5)=0$.
\section{Proof of Theorem~\ref{thm FG'}}\label{sec proof FG'}
We will request the famous Euler's formula~\cite{Andrews}
\begin{equation}\label{Euler-1}
(q;q)_\infty = 1+\sum_{n\geq 1}(-1)^n q^{\fr{3n^2+n}{2}} +\sum_{n\geq 1}(-1)^n q^{\fr{3n^2-n}{2}}
\end{equation}
and the following formula which is found in Fine~\cite[(25.94)]{Fine}
\begin{equation}\label{Fine}
\sum_{n\geq 0}\fr{(aq^{n+1};q)_n t^n}{(q;q)_n} = (t;q)_\infty^{-1} \sum_{n\geq 0}\fr{(t;q)_n}{(q;q)_n} (-at)^n q^{\fr{3n^2+n}{2}}.
\end{equation}
As for part~(a), we have
\[
%\sum_{n\geq 1} F'(n)q^n
\sum_{n\geq 1} q^n (q^{n+1};q)_\infty (q^n;q)_n
=(q;q)_\infty \sum_{n\geq 1} \fr{q^n (q^n;q)_n}{(q;q)_n}
\]
\[
=(q;q)_\infty \sum_{n\geq 0} \fr{q^n (q^n;q)_n}{(q;q)_n} - (q;q)_\infty
\]
\[
=(q;q)_\infty (q;q)_\infty^{-1}\sum_{n\geq 0}(-1)^n q^{\fr{3n^2+n}{2}} - (q;q)_\infty
\]
\[
=\sum_{n\geq 0}(-1)^n q^{\fr{3n^2+n}{2}} - 1-\sum_{n\geq 1}(-1)^n q^{\fr{3n^2+n}{2}} -\sum_{n\geq 1}(-1)^n q^{\fr{3n^2-n}{2}}
\]
\[
=\sum_{n\geq 1}(-1)^{n+1} q^{\fr{3n^2-n}{2}},
\]
where in the fourth identity we applied~\eqref{Fine} with $(a,t)=(q^{-1},q)$
and in the fifth identity we applied~\eqref{Euler-1}.

Regarding part~(b), we have
\[
\sum_{n\geq 1} q^n (q^{n+1};q)_\infty (-q^n;q)_n
=(q;q)_\infty \sum_{n\geq 1} \fr{q^n (-q^n;q)_n}{(q;q)_n}
\]
\[
=(q;q)_\infty \sum_{n\geq 0} \fr{q^n (-q^n;q)_n}{(q;q)_n} - (q;q)_\infty
\]
\[
=(q;q)_\infty (q;q)_\infty^{-1}\sum_{n\geq 0}(-1)^n q^{\fr{3n^2+n}{2}} - (q;q)_\infty
\]
\[
=\sum_{n\geq 0} q^{\fr{3n^2+n}{2}} - 1-\sum_{n\geq 1}(-1)^n q^{\fr{3n^2+n}{2}} -\sum_{n\geq 1}(-1)^n q^{\fr{3n^2-n}{2}}
\]
\[
=2 \sum_{n\geq 0}q^{\fr{3(2n+1)^2 + 2n+1}{2}} - \sum_{n\geq 1} (-1)^n q^{\fr{3n^2-n}{2}}
\]
\[
=2\sum_{n\geq 0}q^{6n^2+7n+2} + \sum_{n\geq 1}(-1)^{n+1} q^{\fr{3n^2-n}{2}},
\]
where in the fourth identity we applied~\eqref{Fine} with $(a,t)=(-q^{-1},q)$ and in the fifth identity we applied~\eqref{Euler-1}.
\section{Proof of Theorem~\ref{thm A'}}\label{sec proof A'}
We first recall the following definition~\cite{Gasper-Rahman}
\[
\pFq{3}{2}{a,b,c}{d,e}{q, z}
=\sum_{n=0}^\infty \fr{(a,b,c;q)_n}{(q,d,e;q)_n}z^n.
\]
We need Jacobi's identity~\cite{Andrews}
\begin{equation}\label{Jacobi}
\sum_{n=0}^\infty (-1)^n (2n+1) q^{\fr{n(n+1)}{2}} = (q;q)_\infty^3,
\end{equation}
the following identity of Gasper and Rahman~\cite[(III.10)]{Gasper-Rahman}
\begin{equation}\label{GR-1}
\pFq{3}{2}{a,b,c}{d,e}{q, \fr{de}{abc}}
=\fr{(b,de/ab,de/bc;q)_\infty}{(d,e,de/abc;q)_\infty}
\pFq{3}{2}{d/b,e/b,de/abc}{de/ab,de/bc}{q, b},
\end{equation}
and the following identity of Andrews and Warnaar~\cite[p. 181]{Andrews-Warnaar 2007}
\begin{equation}\label{AW}
\sum_{n=0}^\infty \fr{(-zq;q^2)_n (-z^{-1}q;q^2)_n q^n}{(-q;q)_{2n+1}}
=\sum_{n=0}^\infty \fr{1-z^{2n+1}}{1-z} z^{-n} q^{n(n+1)}.
\end{equation}
We have
\[
\begin{split}
\mathcal{A}_1(q) &:= \sum_{n=1}^\infty A'(n) q^n \\
&=\sum_{n=1}^\infty q^n (q^{n+1};q)_\infty^3 (q^n;q)_n \\
&= (q;q)_\infty^3 \sum_{n=1}^\infty \fr{q^n (q;q)_{2n-1}}{(q;q)_n^3 (q;q)_{n-1}}.
\end{split}
\]
Hence
\[
\begin{split}
\mathcal{A}_1(q^2)
&=(q^2;q^2)_\infty^3 \sum_{n=1}^\infty \fr{q^{2n} (q^2;q^2)_{2n-1}}{(q^2;q^2)_n^3(q^2;q^2)_{n-1}} \\
 &=(q^2;q^2)_\infty^3 \sum_{n=1}^\infty \fr{q^{2n} (q^2;q^4)_{n}(q^4;q^4)_{n-1}}{(q^2;q^2)_n^3(q^2;q^2)_{n-1}} \\
&= (q^2;q^2)_\infty^3 \sum_{n=1}^\infty \fr{q^{2n} (-q^2;q^2)_{n-1}(q;q^2)_n (-q;q^2)_n}{(q^2;q^2)_n^3} \\
&= (q^2;q^2)_\infty^3 \Big( \fr{1}{2}\pFq{3}{2}{-1,q,-q}{q^2,q^2}{q^2, q^2} -\fr{1}{2} \Big).
\end{split}
\]
Now apply~\eqref{GR-1} with $q\to q^2$, $a=-1$, $b=q$, $c=-q$, $d=e=q^2$ to deduce that
\[
\begin{split}
\pFq{3}{2}{-1,q,-q}{q^2,q^2}{q^2, q^2}
&= \fr{(q,-q^3,-q^2;q^2)_\infty}{(q^2,q^2,q^2;q^2)_\infty}\pFq{3}{2}{q,q,q^2}{-q^3,-q^2}{q^2, q} \\
&=\fr{1}{(1+q) (q^2;q^2)_\infty^3}\pFq{3}{2}{q,q,q^2}{-q^3,-q^2}{q^2, q} \\
&= \fr{1}{(q^2;q^2)_\infty^3} \sum_{n=0}^\infty \fr{(q;q^2)_n^2 q^n}{(-q;q)_{2n+1}} \\
&= \fr{1}{(q^2;q^2)_\infty^3} \sum_{n=0}^\infty (-1)^n q^{n^2+n},
\end{split}
\]
where the last assertion follows by setting $z=-1$ in~\eqref{AW}.
Thus we have proved that
\[
(q^2;q^2)_\infty^3 + 2 \mathcal{A}_1(q^2) = \sum_{n=0}^\infty (-1)^n q^{n^2+n}.
\]
Then by Jacobi's identity~\eqref{Jacobi}
\[
\begin{split}
\mathcal{A}_1(q^2)
&= \fr{1}{2}\Big(\sum_{n=0}^\infty (-1)^n q^{n^2+n} - \sum_{n=0}^\infty (-1)^n (2n+1)q^{n^2+n} \Big) \\
&= - \sum_{n=0}^\infty (-1)^n n q^{n^2+n},
\end{split}
\]
which is the desired formula with $q$ replaced by $q^2$.
\section{Proof of Theorem~\ref{thm A''}}\label{sec proof A''}
We will make an appeal to Gauss' idenity
\begin{equation}\label{Gauss}
\psi(q)=(-q;q)_\infty^2 (q;q)_\infty = \fr{(q^2;q^2)_\infty}{(q;q^2)_\infty} =\sum_{n=0}^\infty q^{\fr{n(n+1)}{2}},
\end{equation}
and the following identity of Gasper and Rahman~\cite[(III.9)]{Gasper-Rahman}
\begin{equation}\label{GR-2}
\pFq{3}{2}{a,b,c}{d,e}{q, \fr{de}{abc}}
=\fr{(e/a,de/bc;q)_\infty}{(e,de/abc;q)_\infty}
\pFq{3}{2}{a,d/b,b/c}{d,de/bc}{q, \fr{e}{a}}.
\end{equation}
We have
\[
\begin{split}
\mathcal{A}_2(q) &:= \sum_{n=1}^\infty A''(n) q^n \\
&=\sum_{n=1}^\infty q^n (-q^{n+1};q)_\infty^2 (q^{n+1};q)_\infty (q^n;q)_n \\
&= (-q;q)_\infty^2(q;q)_\infty \sum_{n=1}^\infty \fr{q^n (q;q)_{2n-1}}{(-q;q)_n^2 (q;q)_{n} (q;q)_{n-1}} \\
&=\psi(q) \sum_{n=1}^\infty \fr{q^n (q;q^2)_n (-q;q)_{n-1}}{(q^2;q^2)_n (-q;q)_n} \\
&=\psi(q)  \Big( \fr{1}{2}\pFq{3}{2}{q^{\fr{1}{2}},-q^{\fr{1}{2}},-1}{-q,-q}{q, q} -\fr{1}{2} \Big)
\end{split}
\]
Hence by replacing $q^2$ with $q$,
\[
\psi(q^2) + 2\mathcal{A}_2(q^2)
=\psi(q^2) \pFq{3}{2}{-q,-1,q}{-q^2,-q^2}{q^2, q^2}.
\]
Now apply~\eqref{GR-2} with $q\to q^2$, $a=-q$, $b=-1$, $c=q$, $d=e=-q^2$ to obtain
\[
\begin{split}
\pFq{3}{2}{-q,-1,q}{-q^2,-q^2}{q^2, q^2}
&= \fr{(q,-q^3;q^2)_\infty}{(-q^2,q^2;q^2)_\infty} \pFq{3}{2}{-q,q^2,-q}{-q^2,-q^3}{q^2, q} \\
&= \fr{(q^2;q^4)_\infty}{(1+q) (-q^2,q^2;q^2)_\infty} \pFq{3}{2}{-q,q^2,-q}{-q^2,-q^3}{q^2, q} \\
&= \fr{1}{(1+q) \psi(q^2)}\ \pFq{3}{2}{q^2,-q,-q}{-q^2,-q^3}{q^2, q} \\
&=\fr{1}{\psi(q^2)} \sum_{n=0}^\infty \fr{(-q;q^2)_n^2 q^n}{(-q^2;q^2)_n (-q;q^2)_{n+1}} \\
&=\fr{1}{\psi(q^2)} \sum_{n=0}^\infty \fr{(-q;q^2)_n^2 q^n}{(-q;q)_{2n+1}} \\
&=\fr{1}{\psi(q^2)} \sum_{n=0}^\infty (2n+1) q^{n(n+1)},
\end{split}
\]
where the last formula follows from~\eqref{AW} with $z=1$.
Thus we have proved that
\[
\psi(q^2)+ 2\mathcal{A}_2(q^2)= \sum_{n=0}^\infty (2n+1) q^{n(n+1)}
\]
and therefore by~\eqref{Gauss},
\[
\mathcal{A}_2(q^2) = \fr{1}{2}\sum_{n=0}^\infty (2n+1) q^{n(n+1)}- \fr{1}{2} \sum_{n=0}^\infty q^{n(n+1)}
=\sum_{n=0}^\infty n q^{n(n+1)},
\]
which is the desired identity with $q^2$ replacing $q$.
\section{Proof of Theorem~\ref{thm B'}}\label{sec proof B'}
We will request Euler's identity~\cite{Andrews}
\begin{equation}\label{Euler}
(-q;q)_\infty = \fr{1}{(q;q^2)_\infty}
\end{equation}
and the $q$-binomial theorem
\begin{equation}\label{q-binomial}
\sum_{n=0}^\infty \fr{(a;q)_n}{(q;q)_n} z^n = \fr{(azq;q)_\infty}{(z;q)_\infty}.
\end{equation}
Now from~\eqref{gen B'},~\eqref{basic-facts}, and~\eqref{Euler}, we find
\[
\begin{split}
\sum_{n=1}^\infty B'(n) q^n
&= \sum_{n=1}^\infty q^{n}(-q^{n+1};q)_\infty^2 (q^{n+1};q)_\infty  (q^{n+1};q)_{n} \\
&=(-q;q)_\infty \sum_{n=1}^\infty q^{n}\fr{(q^{2n+2};q^2)_\infty (q^{n+1};q)_{n}}{(-q;q)_n} \\
&=(-q;q)_\infty \sum_{n=1}^\infty q^{n}\fr{(q^{2n+2};q^2)_\infty (q;q)_{2n}}{(q^2;q^2)_{n}} \\
&= \fr{1}{(q;q^2)_\infty} \sum_{n=1}^\infty q^n (q^{2n+2};q^2)_\infty(q;q^2)_n \\
&= \sum_{n=1}^\infty q^n \fr{(q^{2n+2};q^2)_\infty}{(q^{2n+1};q^2)_\infty} \\
&= \sum_{n=0}^\infty q^n \fr{(q^{2n+2};q^2)_\infty}{(q^{2n+1};q^2)_\infty} - \fr{(q^2;q^2)_\infty}{(q;q^2)_\infty} \\
&= \fr{(q^2;q^2)_\infty}{(q;q^2)_\infty} \sum_{n=0}^\infty \fr{(q;q^2)_n}{(q^2;q^2)_n} q^n
- \fr{(q^2;q^2)_\infty}{(q;q^2)_\infty} \\
&= \Big(\fr{(q^2;q^2)_\infty}{(q;q^2)_\infty} \Big)^2 - \fr{(q^2;q^2)_\infty}{(q;q^2)_\infty} \\
&= \Big(\sum_{n=0}^\infty q^{\fr{n(n+1)}{2}} \Big)^2 - \sum_{n=0}^\infty q^{\fr{n(n+1)}{2}},
\end{split}
\]
where the penultimate identity follows by~\eqref{q-binomial} and the last identity follows from~\eqref{Gauss}.
\section{Proof of Theorem~\ref{thm C'}}\label{sec proof C'}
Recall that a pair of sequences
$(\alpha_n,\beta_n)_{n\geq 0}$ is called a Bailey pair relative to $a$
%(or relative to $(a,q)$ to avoid confusion)
if~\cite{Andrews 1986, Warnaar 2009}
\[
\beta_n = \sum_{r=0}^n \fr{\alpha_r}{(q)_{n-r}(aq)_{n+r}}.
\]
We shall require the following lemma which is an equivalent variant of
Lovejoy~\cite[(1.12)]{Lovejoy 2012}.
\begin{lemma}\label{thm ConjBailey}
If $(\alpha_n,\beta_n)$ is a Bailey pair relative to $a^2$, then
\[
(q;q)_{\infty} (-aq;q)_{\infty}^2 \sum_{n= 0}^\infty \fr{q^n (a;q)_n (a^2 q;q^2)_n}{(-a q;q)_n}\beta_n
\]
\begin{equation}\label{main id}
=  (1-a) \sum_{r,n=0}^\infty \fr{1+aq^{r+2n+1}}{1-aq^r} a^{2n} q^{2n^2+2nr+n+r} \alpha_r.
\end{equation}
\end{lemma}
We want to apply Lemma~\ref{thm ConjBailey} to
the following Bailey pair relative to $q^2$ which can be found in Lovejoy~\cite{Lovejoy 2012}
\begin{equation}\label{BP-D}
\alpha_n = q^{n^2+n}\fr{1-q^{2n+2}}{1-q^2},\quad
\beta_n = \fr{1}{(q;q)_n (q^2; q)_n}.
\end{equation}
Then by Lemma~\ref{thm ConjBailey} with $a=-q$, we find
\[
(1-q)(q;q)_\infty (q^2;q)_\infty^2 \sum_{n=0}^\infty\fr{q^n (-q;q)_n (q^3;q^2)_n}{(q^2;q)_n (q;q)_n (q^2;q)_n}
\]
\begin{equation}\label{BP-D 1}
= \sum_{r,n=0}^\infty
q^{2n^2+2nr+r^2+3n+2r} (1-q^{r+1})(1-q^{2n+r+2}).
\end{equation}
Then by~\eqref{basic-facts} and~\eqref{Euler} the left hand-side of~\eqref{BP-D 1} equals
\[
\sum_{n=0}^\infty q^n (-q;q)_n (q;q^2)_{n+1} (q^{n+1};q)_\infty (q^{n+2};q)_\infty^2
\]
\[
=(-q;q)_\infty \sum_{n=0}^\infty q^n\fr{(q;q^2)_{n+1} (q^{n+1};q)_\infty (q^{n+2};q)_\infty^2}{(-q^{n+1};q)_\infty}
\]
\[
=\fr{1}{(q;q^2)_\infty} \sum_{n=0}^\infty q^n\fr{(q;q^2)_{n+1} (q^{n+1};q)_\infty (q^{n+2};q)_\infty^2}{(-q^{n+1};q)_\infty}
=\sum_{n=0}^\infty q^n\fr{(q^{n+1};q)_\infty (q^{n+2};q)_\infty^2}{(-q^{n+1};q)_\infty (q^{2n+3};q)_\infty}
\]
\[
=\sum_{n=0}^\infty q^n\fr{(q^{n+1};q)_\infty^2 (q^{n+2};q)_\infty^2}{(q^{2n+2};q^2)_\infty (q^{2n+3};q^2)_\infty}
=\sum_{n=0}^\infty q^n\fr{(q^{n+1};q)_\infty^2 (q^{n+2};q)_\infty^2}{(q^{2n+2};q)_\infty}
\]
\[
=\sum_{n=0}^\infty q^n (q^{n+1};q)_\infty (q^{n+1};q)_{n+1}(q^{n+2};q)_\infty^2
=\fr{1}{q} \sum_{n=1}^\infty q^n (q^{n};q)_\infty (q^{n};q)_{n}(q^{n+1};q)_\infty^2
\]
\[
=\fr{1}{q} \sum_{n=1}^\infty C'(n) q^n,
\]
where the last identity follows from~\eqref{gen C'}. Then~\eqref{BP-D 1} means that
\[
\sum_{n=1}^\infty C'(n) q^n
= \sum_{r,n=0}^\infty
q^{2n^2+2nr+r^2+3n+2r+1} (1-q^{r+1})(1-q^{2n+r+2})
\]
\[
= \sum_{r,n=0}^\infty (q^{2n^2+2nr+r^2+3n+2r+1} + q^{2n^2+2nr+r^2+5n+4r+4})
\]
\[
- \sum_{r,n=0}^\infty (q^{2n^2+2nr+r^2+5n+3r+3} + q^{2n^2+2nr+r^2+3n+3r+2})
\]
\[
=\sum_{n=0}^\infty q^{2n^2+3n+1}
+ 2 \sum_{r,n=0}^\infty q^{2n^2+2nr+r^2+5n+4r+4}
\]
\[
- \sum_{r,n=0}^\infty (q^{2n^2+2nr+r^2+5n+3r+3} + q^{2n^2+2nr+r^2+3n+3r+2})
\]
\[
=\sum_{n=0}^\infty q^{\fr{1}{8} \big( (4n+3)^2-1\big)}
+2 \sum_{r,n=0}^\infty q^{\fr{1}{8} \big( (2r+3)^2 +(2r +4n+5)^2-2\big)}
\]
\[
- \sum_{r,n=0}^\infty q^{\fr{1}{8} \big( (2r+3)^2 +(2r +4n+3)^2-2\big)}
- \sum_{r,n=0}^\infty q^{\fr{1}{8} \big( (2r+1)^2 +(2r +4n+5)^2 -2\big)}.
\]
Hence
\[
\sum_{n=1}^\infty C'(n) q^{8n+2} = \sum_{n=0}^\infty q^{(4n+3)^2+1}
+2 \sum_{r,n=0}^\infty q^{(2r+3)^2 +(2r +4n+5)^2}
\]
\[
-\sum_{r,n=0}^\infty q^{(2r+3)^2 +(2r +4n+3)^2}
-\sum_{r,n=0}^\infty q^{(2r+1)^2 +(2r +4n+5)^2}
\]
\[
=\sum_{n=0}^\infty q^{(4n+3)^2+1} + 2\sum_{r,n=0}^\infty (-1)^{n+1} q^{(2r+3)^2 + (2r+3+2n)^2}
\]
\[
+\sum_{r,n=0}^\infty q^{(2r+3)^2 + (2r+3+4n)^2} - \sum_{r,n=0}^\infty q^{(2r+1)^2 + (2r+1+4n+4)^2}
\]
\[
=2\sum_{r,n=0}^\infty (-1)^{n+1} q^{(2r+3)^2 + (2r+3+2n)^2}
+\sum_{n=0}^\infty q^{(4n+3)^2+1} - \sum_{n=0}^\infty q^{(4n+1)^2+1} + \sum_{r=1}^\infty q^{2(2r+1)^2}
\]
\[
=2\sum_{r,n=0}^\infty (-1)^{n+1} q^{(2r+3)^2 + (2r+3+2n)^2}
+\sum_{n=0}^\infty (-1)^{n+1}q^{(2n+1)^2+1} + \sum_{r=1}^\infty q^{2(2r+1)^2},
\]
which is the desired formula.
\section{Proof of Theorem~\ref{thm D'}}\label{sec proof D'}
We have the following Bailey pair relative to $1$ which is due to Slater~\cite[H(1) p. 468]{Slater 1951}
\begin{equation}\label{BP-B}
\alpha_n = \begin{cases}
1, & \text{if $n=0$}, \\
q^{n^2} (q^n - q^{-n}), & \text{if $n=0$}
\end{cases},
\qquad
\beta_n = \fr{q^n}{(q;q)_n^2}.
\end{equation}
Then by Lemma~\ref{thm ConjBailey} applied to~\eqref{BP-B} with $a=-1$, we get 
\[
(q;q)_\infty^3\sum_{n=0}^\infty \fr{q^{n} (-1;q)_n (q;q^2)_n}{(q;q)_n}\fr{q^n}{(q;q)_n^2}
\]
\[
=2\sum_{n=0}^\infty \fr{1-q^{2n+1}}{2} q^{2n^2+n} 
+ 2 \sum_{\substack{n=0\\ r=1}}^\infty\fr{1-q^{r+2n+1}}{1+q^r}q^{2n^2+2nr+n+r} q^{r^2}(q^r-q^{-r})
\]
\[
=\sum_{n=0}^\infty (1-q^{2n+1})q^{2n^2+n}
+2\sum_{\substack{n=0\\ r=1}}^\infty \fr{1-q^{r+2n+1}}{1+q^r} q^{2n^2+2nr+r^2+n}(q^{2r}-1),
\]
or equivalently,
\begin{equation}\label{help B-1}
\fr{(q;q)_\infty^3}{2} + (q;q)_\infty^3 \sum_{n=1}^\infty \fr{q^{2n} (-q;q)_{n-1} (q;q^2)_n}{(q;q)_n^3}
\end{equation}
\[
= \fr{1}{2}\sum_{n=0}^\infty (1-q^{2n+1})q^{2n^2+n}
- \sum_{\substack{n=0\\ r=1}}^\infty (1-q^r)(1-q^{r+2n+1}) q^{2n^2+2nr+r^2+n}
\]
\[
= \fr{1}{2}\sum_{n=0}^\infty (1-q^{2n+1})q^{2n^2+n}
- \sum_{\substack{n=0\\ r=1}}^\infty \big( q^{2n^2+2nr+r^2+n} + q^{2n^2+2nr+r^2+3n+2r+1}\big)
\]
\[
+\sum_{\substack{n=0\\ r=1}}^\infty \big( q^{2n^2+2nr+r^2+n+r} + q^{2n^2+2nr+r^2+3n+r+1}\big)
\]
\[
= \fr{1}{2}\sum_{n=0}^\infty (1-q^{2n+1})q^{2n^2+n}
-\sum_{\substack{n=0\\ r=0}}^\infty \big( q^{2n^2+2nr+r^2+n} + q^{2n^2+2nr+r^2+3n+2r+1}\big)
\]
\[
+ \sum_{\substack{n=0\\ r=0}}^\infty \big( q^{2n^2+2nr+r^2+n+r} + q^{2n^2+2nr+r^2+3n+r+1}\big)
\]
\[
= \fr{1}{2}\sum_{n=0}^\infty (1-q^{2n+1})q^{2n^2+n}
- \sum_{n=0}^\infty q^{2n^2 +n} - 2 \sum_{n,r=0}^\infty q^{2n^2+2nr+r^2+3n+2r+1}
\]
\[
+ \sum_{n,r=0}^\infty \big( q^{2n^2+2nr+r^2+n+r} + q^{2n^2+2nr+r^2+3n+r+1}\big)
\]
%\[
%= \sum_{n=0}^\infty q^{\fr{1}{8} ( (4n+1)^2 -1)}
%+ 2 \sum_{n,r=0}^\infty q^{\fr{1}{8} ( (2r+1)^2 + (2r+4n+3)^2-2 )}
%\]
%\[
%- \sum_{n,r=0}^\infty q^{\fr{1}{8}( (2r+1)^2 + (2r+4n+1)^2-2 )}
%- \sum_{n,r=0}^\infty q^{\fr{1}{8}( (2r-1)^2 + (2r+4n+3)^2-2 )}.
%\]
\[
= -\fr{1}{2}\sum_{n=0}^\infty q^{\fr{1}{8}( (4n+1)^2-1 )}
-\fr{1}{2}\sum_{n=0}^\infty q^{\fr{1}{8}( (4n+3)^2-1 )}
-2 \sum_{n,r=0}^\infty q^{\fr{1}{8}( (2r+1)^2 + (2r+4n+3)^2-2 )}
\]
\[
+ \sum_{n,r=0}^\infty q^{\fr{1}{8}( (2r+1)^2 + (2r+4n+1)^2-2 )}
+ \sum_{n,r=0}^\infty q^{\fr{1}{8}( (2r-1)^2 + (2r+4n+3)^2-2 )}
\]
\[
= -\fr{1}{2}\sum_{n=0}^\infty q^{\fr{1}{8}( (4n+1)^2-1 )}
+\fr{1}{2}\sum_{n=0}^\infty q^{\fr{1}{8}( (4n+3)^2-1 )}
- \sum_{n=0}^\infty q^{\fr{1}{8}( 2(2n+1)^2-2 )}
\]
\[
-2 \sum_{n,r=0}^\infty q^{\fr{1}{8}( (2r+1)^2 + (2r+4n+3)^2-2 )}
+2 \sum_{n,r=0}^\infty q^{\fr{1}{8}( (2r+1)^2 + (2r+1+4n)^2-2 )}
\]
%\[
%+ 2 \sum_{n,r=0}^\infty q^{\fr{1}{8}( (2r+1)^2 + (2r+1+4n)^2-2 )}
%- \sum_{n=0}^\infty q^{\fr{1}{8}( 2(2n+1)^2-2 )} +\sum_{n=0}^\infty q^{\fr{1}{8}( (4n+3)^2-1 )}
%\]
which by virtue of~\eqref{Jacobi} means
\begin{equation}\label{help B-2}
(q;q)_\infty^3 \sum_{n=1}^\infty \fr{q^{2n} (-q;q)_{n-1} (q;q^2)_n}{(q;q)_n^3}
\end{equation}
\[
=2 \sum_{n,r=0}^\infty q^{\fr{1}{8}( (2r+1)^2 + (2r+1+4n)^2-2 )}
-2 \sum_{n,r=0}^\infty q^{\fr{1}{8}( (2r+1)^2 + (2r+4n+3)^2-2 )}
\]
\[
-\fr{1}{2}\sum_{n=0}^\infty q^{\fr{1}{8}( (4n+1)^2-1 )}
+\fr{1}{2}\sum_{n=0}^\infty q^{\fr{1}{8}( (4n+3)^2-1 )}
- \sum_{n=0}^\infty q^{\fr{1}{8}( 2(2n+1)^2-2 )}
\]
\[
-\fr{1}{2}\sum_{n=0}^\infty (-1)^n (2n+1) q^{\fr{n(n+1)}{2}}
\]
\[
=2 \sum_{n,r=0}^\infty (-1)^nq^{\fr{1}{8}( (2r+1)^2 + (2r+1+2n)^2-2 )}
-\fr{1}{2}\sum_{n=0}^\infty (-1)^n q^{\fr{1}{8}( (2n+1)^2-1 )}
\]
\[
-\sum_{n=0}^\infty q^{\fr{1}{8}( 2(2n+1)^2-2 )}
-\fr{1}{2}\sum_{n=0}^\infty (-1)^n (2n+1) q^{\fr{n(n+1)}{2}}.
\]
Furthermore, by~\eqref{basic-facts} and~\eqref{gen D'}, the left hand-side of~\eqref{help B-2} equals
\[
\sum_{n=1}^\infty q^{2n}\fr{(-q;q)_{n-1} (q;q)_{2n} (q^{n+1};q)_\infty^3}{(q^2;q^2)_n}
= \sum_{n=1}^\infty q^{2n}\fr{(q^2;q^2)_{n-1} (q^n;q)_{n+1} (q^{n+1};q)_\infty^3}{(q^2;q^2)_n}
\]
\begin{equation}\label{help B-3}
= \sum_{n=1}^\infty q^{2n} (q^n;q)_n (q^{n+1};q)_\infty^3 =  \sum_{n=1}^\infty D'(n) q^n.
\end{equation}
Now combining~\eqref{help B-2} and~\eqref{help B-3}, we achieve
\[
\sum_{n=1}^\infty D'(n) q^{8n+2}
= 2 \sum_{n,r=0}^\infty (-1)^n q^{(2r+1)^2 + (2r+1+2n)^2}
-\fr{1}{2}\sum_{n=0}^\infty (-1)^n q^{(2n+1)^2+1}
\]
\[
-\sum_{n=0}^\infty q^{2(2n+1)^2} - \fr{1}{2} \sum_{n=0}^\infty (-1)^n (2n+1) q^{(2n+1)^2 +1}
\]
\[
= 2 \sum_{n,r=0}^\infty (-1)^n q^{(2r+1)^2 + (2r+1+2n)^2}
-\sum_{n=0}^\infty (-1)^n(n+1) q^{(2n+1)^2 +1} - \sum_{n=0}^\infty q^{2(2n+1)^2},
\]
which is the desired formula.
\section{Concluding remarks}\label{sec conclusion}
{\bf 1.\ }Our proofs for the identities in Theorems~\ref{thm FG'}-\ref{thm D'} rely on the theory of $q$-series. It would be very interesting to find combinatorial proofs for these identities.

{\bf 2.\ }
We start with some corollaries of our main results.
We get from Theorems~\ref{thm A'}~and~\ref{thm A''}.
\begin{corollary}\label{cor A'-A''}
For all positive integer $n$ which is not a triangular number, we have
$A'(n)=A''(n)=0$.
\end{corollary}
The following result is a direct consequence of Theorem~\ref{thm B'}.
\begin{corollary}\label{cor B'}
For all positive integer $n$ which is not the sum of two triangular numbers, we have
$B'(n)=0$.
\end{corollary}
Furthermore, noting the elementary fact that an integer $n$ is the sum of two triangular numbers if and only if $8n+2$ is
the sum of two squares, we get the following consequence of Theorem~\ref{thm C'}.
\begin{corollary}\label{cor C'}
For any positive which is not the sum of two triangular numbers, we have
$C'(n)=0$.
\end{corollary}
By the note just before Corollary~\ref{cor C'},  we get the following consequence of Theorem~\ref{thm D'}.
\begin{corollary}\label{cor D'}
For all positive integer $n$ which is not the sum of two triangular numbers, we have
$D'(n)=0$.
\end{corollary}
The vanishing of the differences in Corollary~\ref{cor A'-A''} obviously means that
$A_0(n)=A_1(n)$ and $A_2(n)=A_3(n)$ for all non-triangular numbers and
vanishing of the differences in Corollary in Corollaries~\ref{cor B'},~\ref{cor C'},
and~\ref{cor D'} respectively means that $B_0(n)=B_1(n)$, $C_0(n)=C_1(n)$, and $D_0(n)=D_1(n)$
for any positive integer which is not the sum of two triangular number.
It is natural to ask for bijective proofs for these equalities.

\bigskip
\noindent{\bf Acknowledgment.} The authors are grateful to the referee for valuable comments and interesting suggestions 
which have improved the presentation and quality of the paper.

\noindent{\bf Data Availability Statement.\ }
Not applicable.
\end{document}